\documentclass[11pt]{amsart}

\usepackage{latexsym}
\usepackage{amscd}
\usepackage[all]{xy}
\usepackage[mathscr]{euscript}
\normalsize

\RequirePackage{xspace}

\newcommand{\thought}[1]{}
\renewcommand{\thought}[1]{ \textbf{[#1]}}

\usepackage{enumerate}        
\newenvironment{roenumerate}{\begin{enumerate}[\upshape (i)]}{\end{enumerate}}

\newcommand\nc {\newcommand}
\newcommand\rnc{\renewcommand}

\newcount\blopone
\newcount\xone
\newcount\xtwo
\newcount\ytwo

\newtheorem{theorem}{Theorem}[section]
\newtheorem{prop}[theorem]{Proposition}
\newtheorem{conj}[theorem]{Conjecture}
\newtheorem{com}[theorem]{Comment}
\newtheorem{redu}[theorem]{Reduction}
\newtheorem{refinement}[theorem]{Refinement}
\newtheorem{summary}[theorem]{Summary}
\newtheorem{importnota}[theorem]{Important Notation}
\newtheorem{prblm}[theorem]{Problem}
\newtheorem{notation}[theorem]{Notation}
\newtheorem{defin}[theorem]{Definition}
\newtheorem{caution}[theorem]{Caution}
\newtheorem{remark}[theorem]{Remark}
\newtheorem{reminder}[theorem]{Reminder}
\newtheorem{illustration}[theorem]{Illustration}
\newtheorem{lemma}[theorem]{Lemma}
\newtheorem{construction}[theorem]{Construction}
\newtheorem{corollary}[theorem]{Corollary}
\newtheorem{example}[theorem]{Example}
\newtheorem{conclusion}[theorem]{Conclusion}
\newtheorem{triviality}[theorem]{Triviality}
\newtheorem{proto}[theorem]{Prototype Quasifibration}
\newtheorem{cauex}[theorem]{Cautionary Example}
\newtheorem{hypo}[theorem]{Hypothesis}
\newtheorem{subth}{ }[theorem]
\newtheorem{case}{Case}[theorem]
\newtheorem{ssubth}{ }[subth]
\newtheorem{facts}[theorem]{Facts}

\nc\tri[1]{\begin{triviality}
\label{#1}}
\nc\fac[1]{\begin{facts}
\label{#1}
\begin{em}}
\nc\cas[1]{\begin{case}
\label{#1}
\begin{em}}
\nc\rfn[1]{\begin{refinement}
\label{#1}}
\nc\prt[1]{\begin{proto}
\label{#1}}
\nc\lem[1]{\begin{lemma}
\label{#1}}
\nc\pro[1]{\begin{prop}
\label{#1}}
\nc\thm[1]{\begin{theorem}
\label{#1}}
\nc\cnj[1]{\begin{conj}
\label{#1}}
\nc\cor[1]{\begin{corollary}
\label{#1}}
\nc\dfn[1]{\begin{defin}
\label{#1}}
\nc\sthm[1]{\begin{subth}
\label{#1}}
\nc\exm[1]{\begin{example}
\label{#1}
\begin{em}}
\nc\plm[1]{\begin{prblm}
\label{#1}
\begin{em}}
\nc\rmk[1]{\begin{remark}
\label{#1}
\begin{em}}
\nc\rmd[1]{\begin{reminder}
\label{#1}
\begin{em}}
\nc\ntn[1]{\begin{notation}
\label{#1}
\begin{em}}
\nc\smr[1]{\begin{summary}
\label{#1}
\begin{em}}
\nc\cau[1]{\begin{caution}
\label{#1}
\begin{em}}
\nc\hyp[1]{\begin{hypo}
\label{#1}}
\nc\imn[1]{\begin{importnota}
\label{#1}
\begin{em}}
\nc\rdn[1]{\begin{redu}
\label{#1}
\begin{em}}
\nc\cax[1]{\begin{cauex}
\label{#1}
\begin{em}}
\nc\cmt[1]{\begin{com}
\label{#1}
\begin{em}}
\nc\con[1]{\begin{construction}
\label{#1}
\begin{em}}
\nc\ill[1]{\begin{illustration}
\label{#1}
\begin{em}}
\nc\ssthm[1]{\begin{ssubth}
\label{#1}
\begin{em}}
\nc\cnc[1]{\begin{conclusion}
\label{#1}
\begin{em}}

\nc\elem{\end{lemma}}
\nc\erdn{\end{em}\end{redu}}
\nc\erfn{\end{refinement}}
\nc\eprt{\end{proto}}
\nc\ethm{\end{theorem}}
\nc\ecnj{\end{conj}}
\nc\ecor{\end{corollary}}
\nc\edfn{\end{defin}}
\nc\esthm{\end{subth}}
\nc\epro{\end{prop}}
\nc\etri{\end{triviality}}
\nc\eexm{\end{em}
\end{example}}
\nc\ecmt{\end{em}
\end{com}}
\nc\efac{\end{em}
\end{facts}}
\nc\ermk{\end{em}
\end{remark}}
\nc\ermd{\end{em}
\end{reminder}}
\nc\eill{\end{em}
\end{illustration}}
\nc\eplm{\end{em}
\end{prblm}}
\nc\ecas{\end{em}
\end{case}}
\nc\ecau{\end{em}
\end{caution}}
\nc\ecax{\end{em}
\end{cauex}}
\nc\eimn{\end{em}
\end{importnota}}
\nc\entn{\end{em}
\end{notation}}
\nc\econ{\end{em}
\end{construction}}
\nc\esmr{\end{em}
\end{summary}}
\nc\ehyp{
\end{hypo}}
\nc\ecnc{\end{em}
\end{conclusion}}
\nc\essthm{\end{em}
\end{ssubth}}

\nc\sst{\scriptstyle}
\newcommand{\comment}[1]{}
\newcommand{\ri}{\longrightarrow}

\newcommand{\zz}{{\mathbb Z}}

\newcommand{\K}{{\mathbf K}}
\newcommand{\D}{{\mathbf D}}

\nc\op{^{\hbox{\rm\tiny op}}}
\nc\mth{^{\hbox{\rm\tiny th}}}

\nc\script{\mathscr}
\nc\z{\zeta}
\nc\bc{{\mathbb{BC}}}
\nc\Cset{{\mathbb{C}}}
\nc\aA{{\mathbb{A}}}
\nc\ct{{\script T}}
\nc\cf{{\script F}}
\nc\cg{{\script G}}
\nc\ck{{\script K}}
\nc\ch{{\script H}}
\nc\cl{{\script L}}
\nc\cv{{\script V}}
\nc\ccp{{\script P}}
\nc\cw{{\script W}}
\nc\ce{{\script E}}
\nc\cs{{\script S}}
\nc\car{{\script R}}
\nc\cd{{\script D}}
\nc\cc{{\script C}}
\nc\ca{{\script A}}
\nc\ci{{\script I}}
\nc\cj{{\script J}}
\nc\co{{\script O}}
\nc\cx{{\script X}}
\nc\cy{{\script Y}}
\nc\bd{\begin{description}}
\nc\ed{\end{description}}
\nc\ctob{{\script C}at\big(\ci^{op},\ca\big)}
\nc\clim{{\ds\mathop{\rm lim}_{\ds\longleftarrow}}\,}
\nc\climi{\clim^{\!i}\,}
\nc\climn{\clim^{\!n}\,}
\nc\colim{{\ds\mathop{\rm colim}_{\ds\la}}}
\nc\oa{\overline{\ca}}
\nc\s{\sigma}
\nc\ta{\tau}
\nc\os{\overline\sigma}
\nc\ot{\overline\tau}
\nc\T{\Sigma}
\nc\Tm{\Sigma^{-1}}
\nc\de[1]{{\mathop{\rm deg(#1)}}}
\nc\Ad[1]{\mathop{\rm Ad}(#1)}
\nc\ad[1]{\mathop{\rm ad}(#1)}
\nc\kth{{\it K}--theory}

\def\der #1 {D\left(#1\right)}
\nc\prf{\begin{proof}}
\nc\eprf{\end{proof}}
\nc\ds{\displaystyle}
\nc\Tor{\text{\rm Tor}}
\nc\cb{{\script B}}
\nc\ab{{\script A}b}
\nc\be{\begin{roenumerate}}
\nc\ee{\end{roenumerate}}
\nc\cat[1]{{\script C}at\Big({\big\{#1\big\}}\op\,\,,\,\,\ab\Big)}
\nc\csab{{\script C}at\big(\cs^{op},\ab\big)}
\nc\ctab{{\script C}at\Big({\{\ct^\alpha\}}^{op},\ab\Big)}
\nc\csex{{\script E}x\big(\cs^{op},\ab\big)}
\nc\ctex{{\script E}x\Big({\{\ct^\alpha\}}^{op},\ab\Big)}
\nc\sub{\qquad\subset\qquad}
\nc\ctr[1]{{\left.\ct\left(-,#1\right)\right|}_{\cs}}
\nc\ctrf[2]{{\left.\ct\left(#1,#2\right)\right|}_{\cs}}
\nc\Ctr[1]{{\left.\ct\left(-,#1\right)\right|}_{\ct^\alpha}}
\nc\Ctrf[2]{{\left.\ct\left(#1,#2\right)\right|}_{\ct^\alpha}}
\nc\la{\longrightarrow}
\nc\nin{\noindent}
\nc\cad[1]{\text{card}(#1)}
\nc\eq{\quad=\quad}
\nc\BA{\begin{array}{c}}
\nc\EA{\end{array}}
\nc\barr{
\[
\begin{array}{cccccccccccccccc}
}
\nc\earr{
\end{array}
\]
}
\nc\as[1]{{\langle S\rangle}^{#1}}
\nc\sh{\text{\it shift}}
\nc\yy[1]{{\left.\ct\left(-,#1\right)\right|}_{\ct^c}}
\nc\vrep[2]{{\left.\ct\left(#1,#2\right)\right|}_{\ct^\alpha}}
\nc\da{\downarrow}
\nc\Hom{{\text{\rm Hom}}}
\nc\HHom{{\script H}{\text{\rm om}}}
\nc\RRHom{{\mathbf R}{\script {H}}{\text{\rm om}}}
\nc\RHom{{\mathbf R}{\text{\rm Hom}}}
\nc\Lf{{\mathbf L}}
\nc\R{{\mathbf R}}
\nc\Lt{\,{^\Lf_{}\otimes}\,}
\nc\End{{\mathop{\rm End}}}
\nc\Ext{{\mathop{\rm Ext}}}
\nc\EExt{\script{E}{\mathop{\rm xt}}}
\nc\mr\Modtc
\nc\PExt{{\mathop{\rm PExt}}}
\nc\stm{\text{\rm stmod}(kG)}
\nc\stM{\text{\rm StMod}(kG)}
\nc\e{\varepsilon}
\nc\p{\mathfrak{p}}
\nc\q{\mathfrak{q}}
\nc\rs{\s^{-1}A}
\nc\br{{\{\s^{-1}A\}}}
\nc\ra\ri
\nc\y[1]{\mathbf{y}#1}
\nc\x[1]{\mathbf{z}#1}
\nc\mmod[1]{#1\text{--\rm mod}}
\nc\Mod[1]{#1\text{--\rm Mod}}
\nc\Md {\ensuremath{\mathop{\textup{Mod}}}}
\rnc\mod[1]{\ensuremath{\mathop{\textup{mod-}#1}}\xspace}
\nc\Modtc{\Mod{\ct^c}}
\nc\pgldim[1]{\mathop{\rm pgldim}\,#1}
\nc\tf{{\rm [TR5]}}
\nc\tfs{{\rm [TR5$^*$]}}
\nc\Fun{\text{\rm Funct}(F\op,\ab)}
\nc\sym{\text{\rm Sym}}
\nc\sgn{\text{\rm sgn}}
\nc\Pro{\text{\rm Prod}^{}_\alpha(F\op,\ab)}
\nc\Yt[1]{{\left.\Hom_\ct^{}\left(-,#1\right)\right|}_F^{}}
\nc\dl{\delta}
\nc\Coh{\text{\rm Coh}}
\nc\Fla{\text{\rm Flat}}
\nc\Proj[1]{#1\text{--\rm Proj}}
\nc\proj[1]{#1\text{--\rm proj}}
\nc\Flat[1]{#1\text{--\rm Flat}}
\nc\Inj[1]{#1\text{--\rm Inj}}
\nc\ov{\overline}
\nc\wt{\widetilde}

\nc\ph{\varphi}
\nc\qc{\text{\rm Qcoh}}
\nc\dqc[1]{\D(\qc/#1)}
\nc\dbc[1]{\D^b(\text{\rm Coh}/#1)}
\nc\coh{\text{\rm Coh}}
\nc\spec[1]{\text{\rm Spec}(#1)}
\nc\gen[2]{{\langle#1\rangle}_{#2}}
\nc\bgen[2]{\ov{\langle#1\rangle}_{#2}}
\nc\Coprod{\text{\rm Coprod}}
\nc\Ac{{\mathbf A}}
\nc\hoco{
\begin{picture}(40,10)
\put(20,0){\makebox(0,0)[b]{\text{\rm Hocolim}}}
\put(5,-2){\vector(1,0){30}}
\end{picture}\,\,}

\begin{document}

\author{Carles Casacuberta}\thanks{This article was written during a thematic year 
at the CRM Barcelona. The authors acknowledge support from the
Spanish Ministry of Education and Science under sabbatical grant SAB2006-0135
and research projects MTM2004-03629 and MTM2007-63277. The second-named author was also partly supported
by the Australian Research Council.}
\address{Departament d'\`Algebra i Geometria \\
        Universitat de Barcelona \\
        Gran Via de les Corts Catalanes, 585 \\
        08007 Barcelona \\
        Spain}
\email{Carles.Casacuberta@ub.edu}

\author{Amnon Neeman}\thanks{}
\address{Centre for Mathematics and its Applications \\
        Mathematical Sciences Institute \\
        John Dedman Building \\
        The Australian National University \\
        Canberra, ACT 0200 \\
        Australia}
\email{Amnon.Neeman@anu.edu.au}

\title{Brown representability does not come for free}

\begin{abstract}
We exhibit a triangulated category $\ct$ having both products and
coproducts and a triangulated subcategory $\cs\subset \ct$ which
is both localizing and colocalizing, and for which
neither a Bousfield localization nor a colocalization exists.
It follows that neither the category $\cs$ nor its dual satisfy Brown representability.
Our example involves an abelian category whose derived category
does not have small Hom-sets.
\end{abstract}


\maketitle

\section*{Introduction}
\label{S0}

In recent years, several authors have proved remarkable generalizations
of Brown's representability theorem~\cite{Brown62}; see, for example,
\cite{Franke97,Krause02,Neeman96A,Neeman99}. It therefore becomes important
to have an example of a triangulated category where Brown 
representability fails. In~this short note we produce such a category.

There has also been considerable activity on the subject of
localization in homotopy theory, and in particular on Bousfield's problem of proving 
the existence of localization of spaces with respect to
cohomology theories. In \cite{Casacuberta-Scevenels-Smith05} it was shown
that the existence of cohomological localizations follows from a suitable large-cardinal axiom,
although Bousfield's problem remains open under the ZFC axioms alone.
In a similar vein, it was asked in
\cite[p.\ 35]{Hovey-Palmieri-Strickland} if every localizing subcategory 
(i.e., one which is closed under triangles and coproducts) of a stable homotopy category
admits a Bousfield localization. Although the answer is not known in ZFC either,
the counterexample displayed in the present article shows that Bousfield localizations
need not exist for localizing subcategories of arbitrary triangulated categories.

More explicitly, we show that there is an abelian category $\ca$, due to Freyd, for which
the following holds:

\be
\item
The category $\ca$ satisfies the [AB5] and [AB$4^*$] conditions (it has exact products
and coproducts, and filtered colimits are exact).
\item
Nevertheless, the derived category $\D(\ca)$ does not
have small $\Hom$-sets. That is, there is a proper class of morphisms between 
certain objects of~$\D(\ca)$.
\item
Let $\K(\ca)$ be the homotopy category of chain complexes
in $\ca$, and let $\Ac(\ca)$ be the full subcategory of acyclic
complexes. Then $\Ac(\ca)$ is both a localizing and a
colocalizing subcategory, but neither a Bousfield localization
nor a colocalization exist for $\Ac(\ca)$ in $\K(\ca)$.
\item
Neither the category $\Ac(\ca)$ nor its dual satisfy Brown representability. 
\ee

\section{Description and proof}
\label{S1}

In his 1964 book \cite[Chapter~6, Exercise~A, pp.\ 131--132]{Freyd66},
Freyd constructed an interesting abelian category. We briefly paraphrase
the construction. In this article, our foundational formalism for categories is that of
Mac Lane \cite[I.6]{MacLane71}.

Let $I$ be the class of all small ordinals, and let $R=\zz[I]$ be the 
polynomial ring freely generated by $I$. The ring $R$ has a proper
class of elements, but for what we will do this is no problem.
Let $\ca$ be the abelian category of all small $R$-modules. Thus 
an object in $\ca$ is a small abelian group $M$ together with endomorphisms
$\ph_i\colon M\la M$ for every $i\in I$, such that all the $\ph_i$ commute.
The morphisms in $\ca$ are the $R$-module homomorphisms. Given
two objects $M$ and $N$ in~$\ca$, there is only a small set of
morphisms $\Hom_\ca(M,N)$; it is a subset of the set of abelian group homomorphisms.

Note that the abelian category $\ca$ has many
good properties. It satisfies the [AB5] and [AB$4^*$] conditions.
After all, it is the category of modules over a ring, albeit a
very large ring. However, there is no generator or cogenerator, and 
it will follow from our remarks that there are not enough projectives or injectives.

Let $\zz\in\ca$ be the trivial $R$-module. Thus
the underlying abelian group is the additive group of integers $\zz$, 
and all the maps $\ph_i\colon\zz\la\zz$ are zero.

The following observation is due to Freyd \cite{Freyd66}.

\lem{L1.3} 
With the notation as above, $\Ext^1_\ca(\zz,\zz)$ is a proper class.
\elem

\prf
For every ordinal $i\in I$ we construct a module $M_i$ such that, as an abelian 
group, $M_i=\zz\oplus\zz$. The endomorphisms $\ph_j\colon M_i\la M_i$ are given by
the following rule:
\be
\item 
If $j\neq i$, then $\ph_j\colon M_i\la M_i$ is zero.
\item
The map $\ph_i\colon M_i\la M_i$ is determined by the matrix
\[
\left(
\begin{array}{cc}
0 & 1\\
0 & 0
\end{array}
\right).
\]
\ee
It is clear that the $M_i$ are pairwise non-isomorphic as $R$-modules, since the 
element $j\in I$ for which $\ph_j$ is nonzero on $M_i$ changes as
we change~$i$. Hence, we have a proper class of non-isomorphic modules
$M_i$, each of which fits in an exact sequence
\[
\CD
0 @>>> \zz @>>> M_i @>>> \zz @>>> 0,
\endCD
\]
and we have produced a proper class of elements in $\Ext^1_\ca(\zz,\zz)$.
\eprf

Now consider the category $\K(\ca)$, the homotopy category of~$\ca$. 
The objects are chain complexes of small $R$-modules, and
the morphisms are homotopy equivalence classes of chain maps.
Each $R$-module is viewed as a chain complex concentrated in degree zero.
Let $\Ac(\ca)\subset\K(\ca)$ be the full subcategory of all acyclic
complexes. Both $\K(\ca)$ and
$\Ac(\ca)$ are triangulated categories with small $\Hom$-sets.

In what follows, we refer to \cite{Neeman99} for the necessary terminology and basic facts.
The category $\K(\ca)$ satisfies the [TR5] and [TR$5^*$] conditions;
that is, it has small products and coproducts. The subcategory $\Ac(\ca)$ is \textit{localizing} and
\textit{colocalizing}, meaning that it is closed under both coproducts and products.
(In a triangulated category with coproducts, every triangulated subcategory which is
closed under coproducts is automatically thick by \cite[Proposition~1.6.8]{Neeman99}; 
that is, it contains all direct summands of its objects.) 

The derived category of $\ca$ is the Verdier quotient
\[ 
\D(\ca)=\K(\ca)/\Ac(\ca).
\] 
Since $\Hom_{\D(\ca)}(\zz,\T\zz)\cong\Ext^1_\ca(\zz,\zz)$, Lemma~\ref{L1.3} 
implies the following.

\cor{C1.4}
There is a proper class of morphisms $\zz\la\T\zz$ in $\D(\ca)$.
\qed
\ecor

We remark that it does not help if we restrict attention to bounded derived categories,
since the category $\D^b(\ca)$ does not have small $\Hom$-sets either.

Recall that a \textit{Bousfield localization} for the pair $\Ac(\ca)\subset\K(\ca)$
is a right adjoint of the canonical functor $\K(\ca)\la\D(\ca)$, and a
\textit{Bousfield colocalization} is a left adjoint. 
As shown in \cite[Proposition~9.1.18]{Neeman99},
a Bousfield localization exists for the pair $\Ac(\ca)\subset\K(\ca)$
if and only if the inclusion 
\[ 
i\colon\Ac(\ca)\la\K(\ca) 
\] 
has a right adjoint. Dually, a colocalization
exists if and only if $i$ has a left adjoint.

\cor{C1.6}
There is neither a Bousfield localization
nor a Bousfield colocalization for $\Ac(\ca)$ in $\K(\ca)$.
The inclusion functor $i\colon\Ac(\ca)\la\K(\ca)$ has neither a right 
adjoint nor a left adjoint. 
\ecor

\prf
By \cite[Theorem~9.1.16]{Neeman99}, if a Bousfield
localization existed for $\Ac(\ca)\subset\K(\ca)$,
then the quotient category $\D(\ca)=\K(\ca)/\Ac(\ca)$
would be equivalent to a full subcategory $^\perp\Ac(\ca)\subset\K(\ca)$,
namely the one whose objects are those $X$ such that 
\[
\Hom_{\K(\ca)}(A,X)=0
\]
for all $A\in\Ac(\ca)$.
For this, the category $\D(\ca)$ would have to have small $\Hom$-sets.
Since this is not the case by Corollary~\ref{C1.4},
a Bousfield localization cannot exist. Dually, there can be no Bousfield colocalization.
Therefore, By \cite[Proposition~9.1.18]{Neeman99}, the inclusion of
$\Ac(\ca)$ into $\K(\ca)$ has neither a right adjoint nor a left adjoint. 
\eprf

Let $\ab$ denote the category of (small) abelian groups.
A functor from a triangulated category to $\ab$ is called \textit{homological}
if it takes triangles to long exact sequences.
A triangulated category $\ct$ \textit{satisfies Brown representability}
if it has small coproducts and every homological functor $H\colon\ct^{\rm op}\la\ab$
that takes products to products is representable; that is, there is an object $A$
in $\ct$ such that $H$ is naturally isomorphic to $\Hom_{\ct}(-,A)$. (Note that, 
since products in the 
dual category $\ct^{\rm op}$ are coproducts in $\ct$, our assumption is in fact that
$H$ takes coproducts in $\ct$ to products in~$\ab$.)

\cor{C1.7}
Neither the category $\Ac(\ca)$ nor its dual satisfy Brown representability. 
\ecor

\prf
The category $\Ac(\ca)$ has small products and coproducts, and the
inclusion 
\[ 
i\colon \Ac(\ca)\la\K(\ca) 
\] 
respects both. If Brown
representability held for $\Ac(\ca)$, then the inclusion would have a right adjoint
by \cite[Proposition~9.1.19]{Neeman99}.
If Brown representability held for the dual of $\Ac(\ca)$, then a left
adjoint would have to exist. Corollary~\ref{C1.6} tells us that we have neither. 
\eprf

The failure of Brown representability
for $\Ac(\ca)^{\rm op}$ can be displayed more explicitly, without referring
to results in~\cite{Neeman99}, as follows.
(The argument for $\Ac(\ca)$ is similar.)
The functor $\Hom_{\K(\ca)}^{}(\zz,-)$ is a representable functor
from $\K(\ca)$ to~$\ab$. The composite
\begin{equation}
\label{composite}
\CD
\Ac(\ca) @>i>> \K(\ca) @>\Hom_{\K(\ca)}^{}(\zz,-)>> \ab
\endCD
\end{equation}
is a homological functor taking products
to products, and we assert that it is not representable
by any object of~$\Ac(\ca)$. 

Suppose the contrary. If the composite (\ref{composite}) were representable, then
there would exist a map $\ph\colon\zz\la A$ where $A\in\Ac(\ca)$
and such that all other maps from $\zz$ to
acyclic complexes factor uniquely through $\ph$.
Let us complete $\ph$ to a triangle 
\[
\CD
X @>\alpha>> \zz @>\ph>> A @>>> \T X
\endCD
\]
in $\K(\ca)$. Now any morphism $\zz\la\T\zz$ in $\D(\ca)$ can be realized 
as a pair of maps
\begin{equation}
\label{pair}
\xymatrix@C+20pt@R-20pt{
& Y\ar_{\beta}[dl]\ar[dr] & \\ \zz  &  & \T\zz
}
\end{equation}
where $\beta$ is a quasi-isomorphism. This fits into a triangle
\[
\CD
Y @>\beta>> \zz @>\psi>> B @>>> \T Y
\endCD
\]
with $B\in\Ac(\ca)$. Hence, $\psi$ would factor through the
universal map $\ph\colon \zz\la A$, 
and we discover that the above diagram (\ref{pair}) would be equivalent to a diagram
\[
\xymatrix@C+20pt@R-20pt{
& X\ar_{\alpha}[dl]\ar[dr] & \\
\zz  &                     & \T\zz \,.
}
\]
Thus each morphism $\zz\la\T\zz$ in $\D(\ca)$ would be
represented by some map $X\la \T\zz$ in~$\K(\ca)$, where $X$ is fixed.
Since there is only a small set of such maps, we have
contradicted Corollary~\ref{C1.4}.

\bibliographystyle{amsplain}
\bibliography{Casacuberta-Neeman}

\end{document}